\documentclass[12pt]{article}

\usepackage{graphicx}
\usepackage{amsmath}
\usepackage{amsfonts}
\usepackage{amssymb}
\usepackage[french]{babel}

\newtheorem{theorem}{Th\'eor\`eme}
\newtheorem{lemma}{Lemme}
\newtheorem{proposition}{Proposition}

\newtheorem{definition}{D\'efinition\rm}

\newtheorem{example}{Exemple\/}

\newtheorem{consequence}{Cons\'equence}
\newtheorem{conjecture}{Conjecture}

\newenvironment{proof}[1][D\'emonstration]{\textbf{#1.} }{\ \rule{0.5em}{0.5em}}

\title{Alg\`ebres de Lie r\'esolubles r\'eelles alg\'ebriquement rigides
\footnote{Le troisi\`eme auteur (L. G. V.) remercie la Fundaci\'on Ram\'on Areces qui finance sa bourse pr\'edoctorale.}}

\author{J. M. Ancochea Berm\'udez\dag \footnote{e-mail: ancochea@mat.ucm.es}, R. Campoamor-Stursberg\dag \footnote{e-mail: rutwig@mat.ucm.es}\\ L. Garc\'{\i}a Vergnolle\dag \footnote{e-mail: lucigarcia@mat.ucm.es} et M. Goze\ddag \footnote{e-mail: M.Goze@uha.fr}\\
\\
\dag\ Dpto. Geometr\'{\i}a y Topolog\'{\i}a,\\ Facultad CC. Matem\'aticas U.C.M.\\Plaza de Ciencias 3, E-28040 Madrid\\
\ddag Laboratoire de Math\'ematiques et Applications,\\
Universit\'e de Haute Alsace\\4 rue des Fr\`eres Lumi\`ere,
F-68093 Mulhouse cedex }
\date{}

\begin{document}

\maketitle

\begin{abstract}
Nous pr\'esentons toutes les alg\`ebres de Lie r\'eelles,
r\'esolubles et alg\'ebriquement rigides de dimension inf\'erieure
ou \'egale \`a 8. Nous soulignerons les diff\'erences qui
distinguent cette classification de celle des alg\`ebres complexes
 r\'esolubles rigides.
\end{abstract}

\medskip
Mots clefs:
alg\`ebre de Lie, rigide, r\'esoluble.

\section{D\'efinitions et propri\'et\'es pr\'eliminaires}

Soit un corps commutatif $\mathbb{K}$ de caract\'eristique nulle
et $L_{n}$ la vari\'et\'e alg\'ebrique des $\mathbb{K}$-alg\`ebres
de Lie de dimension $n$. On dit qu'une $\mathbb{K}$-alg\`ebre de
Lie $\frak{g}$ est rigide si son orbite dans $L_{n}$, sous
l'action du groupe $GL(n,\mathbb{K})$, est ouverte. Selon le
crit\`ere de rigidit\'e de Nijenhuis et Richardson \cite{NR},
l'annulation du deuxi\`eme groupe de cohomologie
$H^{2}(\frak{g},\frak{g})$ entra\^ine la rigidit\'e. La
r\'eciproque ne se v\'erifie pas en g\'en\'eral et il existe des
contre-exemples pour $\mathbb{K}=\mathbb{C},\mathbb{R}$ et $n>11$.
Par contre toute alg\`ebre de Lie rigide sur $\mathbb{C}$ de
dimension $n\leq 8$
 v\'erifie $H^{2}(\frak{g},\frak{g})=0$.

Soit $\frak{g}$ un alg\`ebre de Lie r\'eelle. Elle est appel\'ee
forme r\'eelle d'une alg\`ebre de Lie complexe $\frak{g}^{\prime}$
si $\frak{g}^{\prime}$ est isomorphe sur $\mathbb{C}$ \`a
l'alg\`ebre $\frak{g}\otimes_{\mathbb{R}}\mathbb{C}$. Dans ce cas
il existe une base de $\frak{g}^{\prime}$ par rapport \`a laquelle
les constantes de structure sont r\'eelles. Si
$\left\{Y_{1},..,Y_{n}\right\}$ est une telle base de
$\frak{g}^{\prime}$, alors
$\left[Y_{i},Y_{j}\right]=C_{ij}^{k}X_{k}$ avec
$C_{ij}^{k}\in\mathbb{R}$. L'alg\`ebre de Lie r\'eelle d\'efinie
par les m\^{e}mes constantes de structure est une forme r\'eelle
de $\frak{g}^{\prime}$.

\begin{definition}
Soit $\frak{g}$ une alg\`ebre de Lie r\'eelle de dimension $n$. On
appelle tore ext\'erieur de d\'erivations toute sous-alg\`ebre
ab\'elienne dont les \'el\'ements sont semi-simples.
\end{definition}

Ceci signifie que les endomorphismes complexes $f\otimes {\rm
Id}\in End(\frak{g}^{\mathbb{C}})$ commutent et sont simultanement
diagonalisables sur $\mathbb{C}$. On en d\'eduit que tous les
tores ext\'erieurs  maximaux de $\frak{g}$ ont la  m\^{e}me
dimension. Cette dimension est le rang de $\frak{g}$. \newline
Dans le cas complexe, le th\'eor\`eme de Mal'cev \cite{Ma}
pr\'ecise que tous les tores sont conjug\'es par un automorphisme
interne de l'alg\`ebre de Lie. Sur le corps r\'eel ceci n'est plus
le cas.

\begin{example}
Soit $\frak{g}=\frak{h}_{1}$ l'alg\`ebre de Heisenberg de dimension 3. Elle est d\'efinie par $\left[X_{1},X_{2}\right]=X_{3}$. Soit $f\in Der(\frak{g})$. Alors la matrice de $f$ sur cette base s'\'ecrit:
$$\left(
\begin{array}{ccc}
a_{1}&b_{1}&0\\
a_{2}&b_{2}&0\\
a_{3}&b_{3}&a_{1}+a_{2}
\end{array}
\right)$$
et $f$ est $\mathbb{C}$-diagonalisable si et seulement si le mineur
$$\left(
\begin{array}{cc}
a_{1}&b_{1}\\
a_{2}&b_{2}
\end{array}
\right)$$
est diagonalisable. Ainsi le tore complexe est engendr\'e par
\begin{equation*}
\left(
\begin{array}{ccc}
\lambda_{1} & 0 & 0\\
0 & 0 &0  \\
0& 0  &\lambda_{1}
\end{array}
\right) \text{ et } \left(
\begin{array}{ccc}
0 & 0&0\\
0& \lambda_{2}&0\\
0&0  &\lambda_{2}
\end{array}
\right) ,
\end{equation*}
\end{example}
Dans le cas r\'eel on a les tores suivantes:

\begin{eqnarray}
\frak{t}_{1}=\mathbb{R}\left\{\left(
\begin{array}{ccc}
\lambda_{1} & 0 &0\\
0 & 0 & 0 \\
0& 0 & \lambda_{1}
\end{array}
\right), \quad \left(
\begin{array}{ccccc}
0 & 0&0\\
0& \lambda_{2}&0\\
0& 0 & \lambda_{2}
\end{array}
\right)\right\} ,\nonumber \\
\frak{t}_{2}=\mathbb{R}\left\{\left(
\begin{array}{ccc}
a & 0 &0\\
0 & a & 0 \\
0& 0 & 2a
\end{array}
\right), \quad \left(
\begin{array}{ccccc}
0 & -b&0\\
b& 0&0\\
0& 0 & 0
\end{array}
\right)\right\},
\end{eqnarray}
et ces deux tores ne sont pas conjug\'es par rapport \`a
$Aut(\frak{h}_{1})$. Notons que si $\frak{g}$ est r\'eelle de rang
non nul le nombre de classes d'automorphisme des tores est fini et
c'est un invariant de $\frak{g}$. On appellera ce nombre {\it
l'indice toro\"{\i}dal} de $\frak{g}$. Ici il vaut 2, et plus
g\'en\'eralement pour l'alg\`ebre de Heisenberg r\'eelle
$\frak{h}_{p}$ il vaut $p+1$.

\medskip

Notons enfin que si l'alg\`ebre de Lie complexe
$\frak{g}^{\prime}$ admet une forme r\'eelle $\frak{g}$, alors
\begin{equation*}
\dim_{\mathbb{R}}H^{2}(\frak{g},\frak{g})=\dim_{\mathbb{C}}H^{2}(\frak{g}^{\prime},
\frak{g}^{\prime}).
\end{equation*}
En effet, si $\left\{X_{1},..,X_{n}\right\}$ est une base de
$\frak{g}$, les constantes de structure de $\frak{g}^{\prime}$ sur
la base $\left\{Y_{i}:=X_{i}\otimes 1, 1\leq i\leq n\right\}$ sont
\'egales \`a ceux de $\frak{g}$. Il s'en suit que le syst\`eme
lin\'eaire d\'efinissant les cocycles (et les cobords) de
$\frak{g}$ et $\frak{g}^{\prime}$ ont le m\^{e}me rang.

\begin{definition}
Une alg\`ebre de Lie r\'eelle $\frak{g}$ est dite alg\'ebriquement
rigide si \[ \dim_{\mathbb{R}}H^{2}(\frak{g},\frak{g})=0.
\]
\end{definition}

On en d\'eduit que $\frak{g}^{\prime}=\frak{g}\otimes\mathbb{C}$ est alg\'ebriquement rigide.

\begin{consequence}
Si $\dim\frak{g}\leq 8$, alors $\frak{g}$ est rigide si et seulement elle est alg\'ebriquement rigide.
\end{consequence}

Le probl\`eme qui se pose alors est de savoir s'il existe des formes r\'eelles $\frak{g}_{1}$ et $\frak{g}_{2}$ non
 isomorphes d'une alg\`ebre de Lie rigide complexe qui soient toutes les deux rigides. L'exemple suivant montre que
 cette situation appara\^{\i}t d\`es la dimension 4.

\section{Exemples}

\subsection{Dimension 4}

Soit $\frak{a}_{2}$ l'alg\`ebre de Lie ab\'elienne r\'eelle de
dimension 2. Elle admet deux tores non conjug\'es:

\begin{eqnarray}
\frak{t}_{1}=\mathbb{R}\left\{\left(
\begin{array}{cc}
\lambda_{1} & 0\\
0 & 0
\end{array}
\right), \quad \left(
\begin{array}{cc}
0 & 0\\
0& 1
\end{array}
\right)\right\} ,\\
\frak{t}_{2}=\mathbb{R}\left\{\left(
\begin{array}{cc}
1 & 0 \\
0 & 1
\end{array}
\right), \quad \left(
\begin{array}{cc}
0 & -1\\
1& 0
\end{array}
\right)\right\}.
\end{eqnarray}
Soit $\frak{g}_{1}$ l'alg\`ebre de Lie r\'esoluble d\'efinie par
l'extension de $\frak{a}_{2}$ par $\frak{t}_{1}$:
$\frak{g}_{1}=\frak{t}_{1}\oplus\frak{a}_{2}$. Les constantes de
structure sont donn\'ees par
\begin{equation}
\begin{array}{llll}
\left[ X_{1},Y_{1}\right] =Y_{1}, & [X_{1},Y_{2}]=0, &
[X_{2},Y_{1}]=0, & \left[ X_{2},Y_{2}\right] =Y_{2}.
\end{array}
\end{equation}
Elle est isomorphe \`a l'alg\`ebre $\frak{r}_{2}^{2}$
\cite{AG}.\newline Soit $\frak{g}_{2}$ l'extension de
$\frak{a}_{2}$ par $\frak{t}_{2}$,
$\frak{g}_{2}=\frak{t}_{2}\oplus\frak{a}_{2}$. Dans ce cas, les
constantes de structure sont donn\'ees par
\begin{equation}
\begin{array}{llll}
\left[ X_{1},Y_{1}\right] =Y_{1}, & [X_{1},Y_{2}]=Y_{2}, &
[X_{2},Y_{1}]=Y_{2}, & \left[ X_{2},Y_{2}\right] =-Y_{1}.
\end{array}
\end{equation}
Les deux alg\`ebres sont alg\'ebriquement rigides et non isomorphes. Elles admettent comme alg\`ebre complexifi\'ee
$\frak{r}_{2}^{2}$.

\subsection{Dimension 5}

Soit $\mathcal{N}_{5,3}$ l'alg\`ebre de Lie nilpotente r\'eelle de dimension 5 d\'efinie par:
\[
[Y_{1},Y_{2}]=Y_{3},\quad [Y_{1},Y_{3}]=Y_{4},\quad [Y_{2},Y_{3}]=Y_{5}.
\]
Par rapport \`a cette base, les d\'erivations ext\'erieures s'\'ecrivent comme
\begin{equation}
\left(
\begin{array}{ccccc}
f_{1}^{1}&f_{2}^{1}&0&0&0\\
f_{1}^{2}&f_{2}^{2}&0&0&0\\
0&0&f_{1}^{1}+f_{2}^{2}&0&0\\
0&f_{2}^{4}&0&2f_{1}^{1}+f_{2}^{2}&f_{2}^{1}\\
f_{1}^{5}&f_{2}^{5}&0&f_{1}^{2}&f_{1}^{1}+2f_{2}^{2}
\end{array}
\right).
\end{equation}
Le tore de $\mathcal{N}_{5,3}\otimes\mathbb{C}$ est de dimension 2
et l'indice toro\"{\i}dal de $\mathcal{N}_{5,3}$ est \'egal \`a 2.
En effet, elle admet les deux tores suivantes qui ne sont pas
conjug\'es:

\begin{enumerate}

\item le tore $\frak{t}_{1}$ est engendr\'e par deux d\'erivations
diagonales:
\begin{eqnarray}
\frak{t}_{1}=\mathbb{R}\left\{\left(
\begin{array}{ccccc}
1 & 0 &0&0&0\\
0&0&0&0&0\\
0&0&1&0&0\\
0&0&0&2&0\\
0 & 0&0&0&1
\end{array}
\right), \quad \left(
\begin{array}{ccccc}
0 & 0&0&0&0\\
0&1&0&0&0\\
0&0&1&0&0\\
0&0&0&1&0\\
0& 0&0&0&2
\end{array}
\right)\right\}
\end{eqnarray}

\item Le tore $\frak{t}_{2}$ est engendr\'e par
\begin{eqnarray}
\frak{t}_{2}=\mathbb{R}\left\{\left(
\begin{array}{ccccc}
1 & 0 &0&0&0\\
0&1&0&0&0\\
0&0&2&0&0\\
0&0&0&3&0\\
0 & 0&0&0&3
\end{array}
\right), \quad \left(
\begin{array}{ccccc}
0 & 1&0&0&0\\
-1&0&0&0&0\\
0&0&0&0&0\\
0&0&0&0&1\\
0& 0&0&-1&0
\end{array}
\right)\right\}
\end{eqnarray}
\end{enumerate}

Comme les matrices $\left(\begin{array}{cc}
f_{1}^{1} & f_{2}^{1}\\
f_{1}^{2} & f_{2}^{2}
\end{array}\right)$ et  $\left(\begin{array}{cc}
2f_{1}^{1}+f_{2}^{2} & f_{2}^{1}\\
f_{1}^{2} & f_{1}^{1}+2f_{2}^{2}
\end{array}\right)$ ont le m\^{e}me polyn\^{o}me caract\'eristique, elles admettent
les m\^{e}mes valeurs propres. On en d\'eduit que $\frak{t}_{1}$
et $\frak{t}_{2}$ sont les seuls tores \`a conjugaison pr\`es.
Consid\'erons les deux alg\`ebres r\'eelles de dimension 7
\begin{eqnarray*}
\frak{g}_{1}=\frak{t}_{1}\otimes\mathcal{N}_{5,3},\\
\frak{g}_{2}=\frak{t}_{2}\otimes\mathcal{N}_{5,3}.
\end{eqnarray*}
Elles sont alg\'ebriquement rigides et non isomorphes dans
$\mathbb{R}$. Les constantes de structure de $\frak{g}_{2}$ sont
donn\'ees par:
\begin{equation}
\begin{array}{llll}
\left[ X_{1},X_{2}\right] =X_{3}, & [X_{1},X_{3}]=X_{4}, &
[X_{2},X_{3}]=X_{5}, & \left[ Y,X_{1}\right] =-X_{2}, \\
\left[ Y,X_{2}\right] =X_{1}, & [Y,X_{4}]=-X_{5}, & [Y,X_{5}]=X_{4}, & \left[Z,X_{1}\right] =X_{1}, \\
\left[ Z,X_{2}\right] =X_{2}, & [Z,X_{3}]=2X_{3}, & [Z,X_{4}]=3X_{4}, &
[Z,X_{5}]=3X_{5}. \\
\end{array}
\end{equation}
Ces deux exemples nous montrent le r\'esultat suivant:

\begin{theorem}
Soit $\frak{n}$ une alg\`ebre de Lie nilpotente de rang non nul
$k$ et d'indice toro\"{\i}dal $p$. Alors, si
$\frak{t}_{1},..,\frak{t}_{p}$ sont les tores non conjug\'es par
rapport au groupe $Aut(\frak{n})$ les alg\`ebres de Lie
$\frak{g}_{i}:=\frak{t}_{i}\oplus\frak{n}$ sont non isomorphes. Si
$\frak{g}^{\mathbb{C}}:=\frak{g}_{i}\oplus\mathbb{C}$ est
alg\`ebriquement rigide, les alg\`ebres r\'eelles $\frak{g}_{i}$
sont alg\`ebriquement rigides et non isomorphes. Ce sont les
seules, \`a isomorphisme pr\`es, ayant $\frak{g}^{\mathbb{C}}$
comme alg\`ebre complexifi\'ee.
\end{theorem}

\begin{consequence}
Le th\'eor\`eme de d\'ecomposabilit\'e de Carles \cite{Ca} n'est
plus valable pour les alg\`ebres de Lie rigides r\'eelles.
Rappelons que cet r\'esultat dit que toute alg\`ebre de Lie rigide
est la somme semi-directe du nilradical\footnote{C'est-\`a-dire,
le id\'eal nilpotent maximal} et un tore ext\'erieur form\'ee par
des \'el\'ements diagonalisables.
\end{consequence}

\section{La classification r\'eelle jusqu'\`a  dimension 8}

Les exemples ci-dessus montrent que les classifications r\'eelle
et complexe diff\`erent. Nous pouvons pour chacune des alg\`ebres
nilpotentes r\'eelles de dimension $5$ reprendre la preuve
pr\'ecedente. Dans ce paragraphe nous proposons la classification
des alg\`ebres de Lie r\'eelles r\'esolubles rigides de dimension
$n\leq 8$ en se basant  d'une part sur la classification complexe
\cite{AG}, et d'autre part sur la classification des alg\`ebres de
Lie r\'eelles nilpotentes de dimension $m\leq 6$ \cite{Ce,Ve}.

\begin{lemma}
Soit $\frak{g}$ une alg\`{e}bre de Lie r\'{e}soluble r\'eelle
alg\'{e}briquement rigide de dimension 8 dont le tore contient au
moins une d\'{e}rivation non-diagonalisable. Alors le nilradical
$\frak{n}$ est de dimension $m\leq6$.
\end{lemma}

\begin{proof}
Supposons que $\frak{g}$ soit de rang 1, c'est-\`{a}-dire,
$\dim\frak{n}=7$. Alors l'espace compl\'{e}mentaire de $\frak{n}$
dans $\frak{g}$ est de dimension 1, donc engendr\'{e} par un
vecteur $X$ tel que l'op\'erateur adjoint $ad\left( X\right)  $
est diagonalisable sur $\mathbb{C}$. D'apr\`{e}s \cite{A2}, pour
toute alg\`{e}bre de Lie r\'{e}soluble complexe rigide, il existe
une base du tore ext\'{e}rieur $\frak{t}$ telle que les valeurs
propres des op\'{e}rateurs adjoints sont enti\`{e}res. On en
d\'{e}duit l'existence d'un vecteur $Y$, tel que $X=\lambda Y$
o\`{u} $\lambda\in\mathbb{C}$ et la forme canonique est
${\rm diag}(n_{1},\dots,n_{7},0)$, o\`{u} $n_{1},\dots,n_{7}%
\in\mathbb{Z}$. De plus, les valeurs propres sont non nulles
d'apr\`es \cite{AG}. L'op\'erateur $ad(X)$ n'\'etant
diagonalisable que sur le corps complexe implique que l'existence
d'une suite de nombres entiers $\left\{n_{1},..,n_{7}\right\}$ et
la partie r\'eelle de $\lambda$ est nulle. Alors $ad\left(
X\right) ={\rm diag}(\lambda n_{1},\dots,\lambda n_{7},0)$ sur
$\mathbb{C}$. Des propri\'et\'es des polyn\^{o}mes \`a
coefficients r\'eels d'ordre impair on d\'eduit que $\lambda
n_{i}$ et $\bar{\lambda}n_{i}$ sont des valeurs propres de
$ad\left(
X\right)$, d'o\`u  $ad(X)={\rm diag}(\lambda n_{1},\bar{\lambda}n_{1},\lambda n_{2}%
,\bar{\lambda}n_{2},\lambda n_{3},\bar{\lambda}n_{3},\lambda n_{4}%
,\bar{\lambda}n_{4})$. Ceci implique qu'il existe au moins un $\lambda
n_{i}=0$, d'o\`{u} la contradiction.
\end{proof}

\smallskip

D'apr\`{e}s le lemme le probl\`{e}me est alors r\'{e}duit \`{a} d\'{e}terminer
les formes r\'{e}elles des alg\`{e}bres de Lie rigides complexes.

\begin{proposition}
Toute alg\`{e}bre de Lie rigide r\'{e}elle de dimension $n\leq8$
possedant au moins une d\'{e}rivation non diagonale est isomorphe
\`{a} une des alg\`{e}bres suivantes:

\begin{itemize}
\item  Dimension 4:%
\[%
\begin{array}
[c]{ccccc}%
\frak{g}_{4}^{2}: & \left[  X_{1},X_{3}\right]  =X_{2}, & \left[  X_{3}%
,X_{2}\right]  =X_{1}, & \left[  X_{4},X_{1}\right]  =X_{1}, & \left[
X_{4},X_{2}\right]  =X_{2}.
\end{array}
\]

\item  Dimension 5:%
\[%
\begin{array}
[c]{ccccc}%
\frak{g}_{5}^{2}: & [X_{1},X_{2}]=X_{3}, & [X_{4},X_{1}]=X_{1}, & [X_{4}%
,X_{2}]=X_{2}, & [X_{4},X_{3}]=2X_{3},\\
& \left[  X_{5},X_{1}\right]  =-X_{2}, & [X_{5},X_{2}]=X_{1}. &  &
\end{array}
\]

\item  Dimension 6:%
\[%
\begin{array}
[c]{ccccc}%
\frak{g}_{6}^{4}: & \left[  X_{1},X_{3}\right]  =X_{2}, & \left[  X_{3}%
,X_{2}\right]  =X_{1}, & \left[  X_{4},X_{1}\right]  =X_{1}, & \left[
X_{4},X_{2}\right]  =X_{2},\\
& \left[  X_{5},X_{6}\right]  =X_{5}. &  &  &
\end{array}
\]

\item  Dimension 7:%
\[%
\begin{array}
[c]{lllll}%
\frak{g}_{7}^{9}: & [X_{1},X_{2}]=X_{3}, & [X_{1},X_{3}]=X_{4}, & [X_{2}%
,X_{3}]=X_{4}\newline , & [X_{6},X_{1}]=-X_{2},\\
& \left[  X_{6},X_{2}\right]  =X_{1}, & \left[  X_{6},X_{4}\right]  =-X_{5}, &
\left[  X_{6},X_{5}\right]  =X_{4}, & \left[  X_{7},X_{1}\right]  =X_{1},\\
& \left[  X_{7},X_{2}\right]  =X_{2}, & \left[  X_{7},X_{3}\right]  =2X_{3}, &
\left[  X_{7},X_{4}\right]  =3X_{4}, & \left[  X_{7},X_{5}\right]  =3X_{5}.
\end{array}
\]

\[%
\begin{array}
[c]{lllll}%
\frak{g}_{7}^{10}: & [X_{5},X_{1}]=X_{1}, & [X_{5},X_{2}]=X_{2}, & [X_{5}%
,X_{3}]=2X_{3}\newline , & [X_{6},X_{1}]=-X_{2},\\
& \left[  X_{6},X_{2}\right]  =X_{1}, & \left[ X_{7},X_{4}\right]
=X_{4}, & & .
\end{array}
\]

\item  Dimension 8:%
\[%
\begin{array}
[c]{ccccc}%
\frak{g}_{8}^{34}: & \left[  X_{1},X_{3}\right]  =X_{2}, & \left[  X_{3}%
,X_{2}\right]  =X_{1}, & \left[  X_{4},X_{1}\right]  =X_{1}, & \left[
X_{4},X_{2}\right]  =X_{2},\\
& \left[  X_{5},X_{7}\right]  =X_{6}, & \left[  X_{7},X_{6}\right]  =X_{5}, &
\left[  X_{8},X_{5}\right]  =X_{5}, & \left[  X_{8},X_{6}\right]  =X_{6}.
\end{array}
\]%
\[%
\begin{array}
[c]{ccccc}%
\frak{g}_{8}^{35}: & \left[  X_{1},X_{3}\right]  =X_{2}, & \left[  X_{3}%
,X_{2}\right]  =X_{1}, & \left[  X_{4},X_{1}\right]  =X_{1}, & \left[
X_{4},X_{2}\right]  =X_{2},\\
& \left[  X_{5},X_{6}\right]  =X_{6}, & \left[  X_{7},X_{8}\right]  =X_{8}. &
&
\end{array}
\]%
\[%
\begin{array}
[c]{ccccc}%
\frak{g}_{8}^{36}: & [X_{1},X_{2}]=X_{4}, & [X_{1},X_{3}]=X_{5}, &
[X_{6},X_{i}]=X_{i}\, & \left(  i=1,4,5\right)  ,\\
& [X_{7},X_{2}]=-X_{3}, & [X_{7},X_{3}]=X_{2}, & [X_{7},X_{4}]=-X_{5}, &
[X_{7},X_{5}]=X_{4},\\
& [X_{8},X_{i}]=X_{i} & \left(  i=2,3,4,5\right)  . &  &
\end{array}
\]%
\[%
\begin{array}
[c]{ccccc}%
\frak{g}_{8}^{37}: & \left[  X_{1},X_{2}\right]  =X_{5}, & \left[  X_{3}%
,X_{4}\right]  =X_{5}, & \left[  X_{6},X_{i}\right]  =X_{i} & \left(
i=1,2\right)  ,\\
& [X_{6},X_{5}]=2X_{5}, & [X_{7},X_{1}]=-X_{2}, & [X_{7},X_{2}]=X_{1}, &
[X_{8},X_{3}]=-X_{4},\\
& [X_{8},X_{4}]=X_{3}\newline . &  &  &
\end{array}
\]%
\[%
\begin{array}
[c]{ccccc}%
\frak{g}_{8}^{38}: & \left[  X_{1},X_{2}\right]  =X_{5}, & \left[  X_{3}%
,X_{4}\right]  =X_{5}, & \left[  X_{6},X_{i}\right]  =X_{i} & \left(
i=1,2\right)  ,\\
& [X_{6},X_{5}]=2X_{5}, & [X_{7},X_{1}]=X_{1}, & [X_{7},X_{2}]=-X_{2}, &
[X_{8},X_{3}]=-X_{4},\\
& [X_{8},X_{4}]=X_{3}\newline . &  &  &
\end{array}
\]%
\[%
\begin{array}
[c]{ccccc}%
\frak{g}_{8}^{39}: & \left[  X_{1},X_{i}\right]  =X_{i+1}, & \left(  2\leq
i\leq4\right)  , & \left[  X_{3},X_{2}\right]  =X_{6}, & \left[  X_{6}%
,X_{2}\right]  =X_{5},\\
& \left[  X_{7},X_{1}\right]  =X_{1}, & \left[  X_{7},X_{2}\right]  =X_{2}, &
\left[  X_{7},X_{3}\right]  =2X_{3}, & \left[  X_{7},X_{4}\right]  =3X_{4},\\
& \left[  X_{7},X_{5}\right]  =4X_{5}, & \left[  X_{7},X_{6}\right]
=3X_{6}, & \left[  X_{8},X_{1}\right]  =X_{2}, & \left[  X_{8},X_{2}\right]
=-X_{1},\\
& \left[  X_{8},X_{4}\right]  =-X_{6}, & \left[ X_{8},X_{6}\right]
=X_{4}. & &
\end{array}
\]%
\[%
\begin{array}
[c]{ccccc}%
\frak{g}_{8}^{40}: & \left[  X_{1},X_{i}\right]  =X_{i+2}, & \left(  2\leq
i\leq4\right)  , & \left[  X_{2},X_{3}\right]  =X_{6}, & \left[  X_{4}%
,X_{2}\right]  =X_{5},\\
& \left[  X_{7},X_{1}\right]  =X_{1}, & \left[  X_{7},X_{2}\right]  =X_{2}, &
\left[  X_{7},X_{3}\right]  =2X_{3}, & \left[  X_{7},X_{4}\right]  =2X_{4},\\
& \left[  X_{7},X_{5}\right]  =3X_{5}, & \left[  X_{7},X_{6}\right]
=3X_{6}, & \left[  X_{8},X_{1}\right]  =X_{2}, & \left[  X_{8},X_{2}\right]
=-X_{1},\\
& \left[  X_{8},X_{5}\right]  =X_{6}, & \left[  X_{8},X_{6}\right]  =-X_{5}. &
&
\end{array}
\]
De plus, les alg\`{e}bres sont deux-\`{a}-deux non isomorphes.
\end{itemize}
\end{proposition}

\begin{theorem}
Soit $\frak{g}$ une alg\`{e}bre de Lie r\'{e}elle, r\'{e}soluble
et alg\'{e}briquement rigide de dimension inf\'{e}rieure ou
\'{e}gale \`{a} 8. Alors $\frak{g}$ est isomorphe \`{a} l'une des
alg\`{e}bres $\frak{g}_{i}^{j}$ d\'{e}crites par les lois
$\mu_{i}^{j}$ de la liste \cite{AG} (celles ci sont les formes
r\'{e}elles obtenues par restriction des scalaires des
alg\`{e}bres complexes r\'{e}solubles  rigides de dimension
inf\'{e}rieure ou \'{e}gale \`{a} 8) ou bien \`{a} l'une des
alg\`{e}bres du lemme pr\'ec\'edent. De plus, toutes ces
alg\`{e}bres sont deux-\`{a}-deux non isomorphes.
\end{theorem}

\begin{proof}
Soit $\frak{g}$ une alg\`{e}bre v\'{e}rifiant les hypoth\`{e}ses
du th\'{e}or\`{e}me.\newline Si la dimension du niradical est
inf\'{e}rieure ou \'{e}gale \`{a} 6, les alg\`ebres poss\'edant
des d\'erivations non diagonalisables sont donn\'ees par le lemme
pr\'ec\'edent, et ceux ayant un tore form\'e par des d\'erivations
diagonalisables sont classifi\'ees dans \cite{AG}. Si la dimension
du niradical est \'{e}gale \`{a} 7, alors n\'{e}cessairement
$\dim\,\frak{g}=8$ et le tore est engendr\'{e} par une seule
d\'{e}rivation $f$. D'apr\`{e}s le lemme ci-dessus, $f$ est
diagonale et donc $\frak{g}$ est obtenue par restriction des
scalaires des alg\`ebres complexes de la liste \cite{AG} qui ont
un nilradical de dimension 7.
\end{proof}

\medskip

De la th\'{e}orie g\'{e}n\'{e}rale de la rigidit\'{e} sur le corps
complexe on d\'eduit que toute alg\`{e}bre de Lie r\'{e}soluble
rigide est determin\'{e}e de fa\c{c}on biunivoque par son
nilradical, et par cons\'{e}quent par le syst\`{e}me de racines
\cite{A2}. Par contre, dans le cas r\'{e}el le nilradical ne
d\'{e}termine pas la structure du tore. La table suivante donne
les alg\`ebres de Lie r\'eelles rigides non isomorphes de
dimension $n\leq 8$ qui s'appuient sur un nilradical donn\'e. On
appelera forme normale \`a l'alg\`ebre de Lie r\'eelle obtenue par
r\'estriction des scalaires de l'alg\`ebre complexe de la
classification dans \cite{AG}.

\begin{table}
\caption{Alg\`ebres de Lie r\'eeles rigides ayant le m\^{e}me
nilradical.\newline }
\begin{tabular}
[c]{|c|c|l|}\hline
Dimension & Nilradical & Formes r\'{e}elles\\\hline
4 & ab\'{e}lien &
\begin{tabular}
[c]{l}%
$\frak{g}_{4}^{1}$\quad (forme normale)\\
$\frak{g}_{4}^{2}$%
\end{tabular}
\\\hline
5 & $\mathcal{N}_{3}\,=\,\frak{h}_{1}$ &
\begin{tabular}
[c]{l}%
$\frak{g}_{5}^{1}$\quad (forme normale)\\
$\frak{g}_{5}^{2}$%
\end{tabular}
\\\hline
6 & ab\'{e}lien &
\begin{tabular}
[c]{l}%
$\frak{g}_{6}^{3}$\quad (forme normale)\\
$\frak{g}_{6}^{4}$%
\end{tabular}
\\\hline
7 & $\mathcal{N}_{3}\oplus\mathbb{R}$ &
\begin{tabular}
[c]{l}%
$\frak{g}_{7}^{8}$\quad (forme normale)\\
$\frak{g}_{7}^{10}$%
\end{tabular}
\\\cline{2-3}%
& $\mathcal{N}_{5,3}$ &
\begin{tabular}
[c]{l}%
$\frak{g}_{7}^{6}$\quad (forme normale)\\
$\frak{g}_{7}^{9}$%
\end{tabular}
\\\hline
 & ab\'{e}lien &
\begin{tabular}
[c]{l}%
$\frak{g}_{8}^{33}\simeq\frak{g}_{2}^{1}\oplus\frak{g}_{2}^{1}\oplus\frak{g}%
_{2}^{1}\oplus\frak{g}_{2}^{1}$\quad (forme normale)\\
$\frak{g}_{8}^{34}$\\
$\frak{g}_{8}^{35}$%
\end{tabular}
\\\cline{2-3}%
& $\mathcal{N}_{5,5}$ &
\begin{tabular}
[c]{l}%
$\frak{g}_{8}^{31}$\quad (forme normale)\\
$\frak{g}_{8}^{36}$%
\end{tabular}
\\\cline{2-3}%
8 & $\mathcal{N}_{5,6}\,=\,\frak{h}_{2}$ &
\begin{tabular}
[c]{l}%
$\frak{g}_{8}^{32}$\quad (forme normale)\\
$\frak{g}_{8}^{37}$\\
$\frak{g}_{8}^{38}$%
\end{tabular}
\\\cline{2-3}%
& $\mathcal{N}_{6,6}$ &
\begin{tabular}
[c]{l}%
$\frak{g}_{8}^{22}$\quad (forme normale)\\
$\frak{g}_{8}^{39}$%
\end{tabular}
\\\cline{2-3}%
& $\mathcal{N}_{6,14}$ &
\begin{tabular}
[c]{l}%
$\frak{g}_{8}^{29}$\quad (forme normale)\\
$\frak{g}_{8}^{40}$%
\end{tabular}
\\\hline
\end{tabular}
\end{table}

\medskip

La raison de l'existence des alg\`ebres rigides r\'eelles
s'appuyant sur le m\^{e}me nilradical est une cons\'equence de
l'invalidit\'e du th\'eor\`eme de d\'ecomposabilit\'e, et plus
pr\'ecisement, de l'existence de d\'erivations ext\'erieures non
diagonalisables. Par cons\'equence, le nombre des alg\`ebres de
Lie r\'eelles rigides ayant le m\^{e}me nilradical est donn\'e
pr\'ecisement par l'indice toro\^{\i}dal d\'efini dans la
premi\`ere section.

\begin{example}
L'alg\`{e}bre de Heisenberg $\frak{h}_{n}$ est d\'{e}finie par les
crochets
\[
\left[X_{2i-1},X_{2i}\right]=X_{2n+1}\quad i=1,\dots,n
\]
Sur le corps complexe, il n'existe qu'une alg\`{e}bre
$\frak{g}_{3n+2}$ rigide ayant $\frak{h}_{n}$ pour nilradical,
donn\'ee par le tore de dimension $n+1$ engendr\'e par $\left\{
Y_{1},..,Y_{n+1}\right\} $:

\begin{equation*}
\begin{array}{lll}
\lbrack Y_{i},X_{2i-1}]=X_{2i-1}, & [Y_{i},X_{2i}]=-X_{2i}, & 1\leq i\leq n
\\
\lbrack Y_{n+1},X_{i}]=X_{i}, & [Y_{n+1},X_{2n+1}]=2X_{2n+1}, &
1\leq i\leq 2n.
\end{array}
\end{equation*}
Les formes r\'{e}elles de $%
\frak{g}_{3n+2}$, que nous d\'{e}notons par $\frak{g}_{3n+2,k}$ (o\`{u} $%
0\leq k\leq n$)  sont donn\'{e}es par:

\begin{equation*}
\begin{array}{lll}
\lbrack X_{2i-1},X_{2i}]=X_{2n+1}, &  & 1\leq i\leq n \\
\lbrack Y_{i},X_{2i-1}]=X_{2i}, & [Y_{i},X_{2i}]=-X_{2i-1}, &
1\leq i\leq
k \\
\lbrack Y_{i},X_{2i-1}]=X_{2i-1}, & [Y_{i},X_{2i}]=-X_{2i}, & k+1\leq i\leq
n+1 \\
\lbrack Y_{n+1},X_{i}]=X_{i}, & [Y_{n+1},X_{2n+1}]=2X_{2n+1} &
1\leq i\leq 2n.
\end{array}
\end{equation*}

En particulier, la forme r\'eelle $\frak{g}_{3n+2,k}$ poss\`ede
$k$ d\'erivations non diagonalisables. \c{C}a montre que l'indice
toro\"{\i}dal est \'egale a $p+1$. Cette \'etude nous conduit \`a
emettre la conjecture suivante:
\end{example}

\begin{conjecture}
Toute alg\`ebre de Lie r\'eelle r\'esoluble rigide poss\`ede au
moins une d\'erivation diagonale.
\end{conjecture}

\section*{Remarques finales}

Dans \cite{C19,C29} on montre que les syst\`emes des poids des
alg\`ebres de Lie nilpotentes \cite{Fa} peuvent \^{e}tre d\'ecrits
par des crit\`eres combinatoires. En particulier, on montre que
toute alg\`ebre de Lie rigide complexe dont l'indice de
r\'esolubilit\'e est deux est d\'ecomposable. L'exemple en
dimension 4 vu dans 2.1 montre que cette propri\'et\'e ne
s'\'etend pas au cas r\'eel. Ceci signifie que la notion de graphe
des poids do\^{\i}t \^{e}tre modifi\'ee pour couvrir la 
classification des rigides r\'eelles. Une autre propri\'et\'e
caract\'eristique du cas complexe concerne la r\'esolubilit\'e
compl\`ete des alg\`ebres rigides. Rappelons qu'une alg\`ebre de
Lie est dite compl\'etement r\'esoluble \cite{Di}  s'il existe une
suite d\'{e}croissante d'id\'{e}aux
\[
\frak{g}=I_{0}\supset I_{1}\supset\ldots I_{n-1}\supset I_{n}=0
\]
telle que $\dim_{\mathbb{K}}I_{k}/I_{k+1}=1$ pour tout $k=0,..,n-1$%
. Sur le corps $\mathbb{C}$, la r\'{e}solubilit\'{e} et la
r\'{e}solubilit\'{e} compl\`{e}te sont equivalentes. Par contre,
sur $\mathbb{R}$, on peut seulement dire que la
r\'{e}solubilit\'{e} compl\`{e}te implique la
r\'{e}solubilit\'{e}, la r\'eciproque \'etant fausse en
g\'en\'eral. D'ailleurs, si $\frak{g}$ une alg\`{e}bre de Lie
r\'{e}soluble rigide complexe alors seulement la forme r\'{e}elle
obtenue par restriction des scalaires (le tenseur de structure
\'etant rationnel) est compl\`{e}tement r\'{e}soluble.


\begin{thebibliography}{99}
\bibitem{A2}J. M. Ancochea, M. Goze, Le rang du syst\`{e}me lin\'{e}aire des
racines d'une alg\`{e}bre de Lie rigide r\'{e}soluble complexe, Comm. Algebra
20 (1992), 875-887.

\bibitem {A3}J. M. Ancochea, M. Goze, On the nonrationality of rigid Lie
algebras, Proc. Am. Math. Soc. \textbf{127} (1999), 2611-2618.

\bibitem {C19}J. M. Ancochea, R. Campoamor-Stursberg. 2-step solvable Lie
algebras and weight graphs, Transf. Groups \textbf{7} (2002), 307-320.

\bibitem {C22}R. Campoamor-Stursberg. Invariants of solvable rigid Lie
algebras up to dimension 8, J. Phys. A: Math. Gen. \textbf{35} (2002), 6293-6306.

\bibitem {C29}R. Campoamor-Stursberg. A graph theoretical determination of
solvable complete rigid Lie algebras, Linear Alg. Appl. \textbf{372} (2003), 53-66.

\bibitem {Ca}R. Carles, Sur la structure des alg\`{e}bres de Lie rigides, Ann.
Inst. Fourier 34 (1984), 65-82.

\bibitem {Ce}A. Cerezo, Les alg\`{e}bres de Lie nilpotentes r\'{e}elles et
complexes de dimension 6, Pr\'{e}publications Universit\'{e} de Nice, 1983.

\bibitem {Di}J. Dixmier. Alg\`{e}bres enveloppantes, Gauthier-Villars, Paris, 1974.

\bibitem{Fa}G. Favre, Syst\`{e}me des poids sur une alg\`{e}bre de Lie
nilpotente, Manuscripta Math. \textbf{9} ( 1973), 53-90.

\bibitem {AG}M. Goze, J. M. Ancochea, On the classification of rigid Lie
algebras, J. Algebra 245 (2001), 68-91.

\bibitem {Ma}A. I. Mal'cev, Solvable Lie algebras, Izv. Akad. Nauk SSSR 9
(1945), 329-356.

\bibitem {NR}A. Nijenhuis, R. W. Richardson, Deformations of Lie algebra
structures. J. Math. Mech. 17 1967 89-105.

\bibitem{Ve} M. Vergne, Vari\'et\'e des alg\`ebres de Lie nilpotentes,
Th\`ese 3\`eme cycle, Paris 1966.
\end{thebibliography}
\end{document}